\documentclass[12pt]{article}
\usepackage{amsfonts, amsmath, amssymb,amsthm,bbm,comment,fancyhdr,float,makeidx,mathrsfs,verbatim}

\normalsize

\newcommand{\si}{\sigma}

\newcommand{\ol}{\overline}

\newcommand{\vs}{\vspace*}

\newcommand{\nin}{\noindent}

\newtheorem{mthm}{Theorem}[section]
\newtheorem{mylem}[mthm]{Lemma}
\newtheorem{myprn}[mthm]{Proposition}
\newtheorem{mycor}[mthm]{Corollary}
\newtheorem{mydef}[mthm]{Definition}
\newtheorem{myrem}[mthm]{Remark}
\newtheorem{mycon}[mthm]{Construction}
\newtheorem{myeg} [mthm]{Example}
\newtheorem{myque} [mthm]{Question}

\def \nin {\noindent}

\def \Lemma #1 {\vs{3mm}\nin {\bf Lemma #1} \it}
\def \Prop #1 {\vs{3mm}\nin {\bf Proposition #1} \it}
\def \Th #1 {\vs{3mm}\nin {\bf Theorem #1} \it}
\def \Cor #1 {\vs{3mm}\nin {\bf Corollary #1} \it}
\def \Ex #1 {\vs{3mm}\nin {\bf Example #1} \it}

\def \part #1 {\hfil\break\hglue 12pt {\rm (#1)~}}

\def \qed {~\vrule height6pt width 6pt depth 0pt}
\def\fs{\footnotesize}

\def\longra{\longrightarrow}

\title{
\bf\LARGE  Minimal free resolution of a graded ideal with linear quotients\thanks{This research is supported by the National Natural
Science Foundation of China (Grant No. 11271250). } }
\author{{A-Ming Liu\thanks{aming8809@163.com} and Tongsuo Wu\thanks{Corresponding author. tswu@sjtu.edu.cn}}\\
 {\small Department of Mathematics, Shanghai Jiaotong University}
}

\date{}
\begin{document}
\baselineskip=16pt \maketitle

\begin{center}
\begin{minipage}{12cm}

 \vs{3mm}\nin{\small\bf Abstract.} {\fs  Let $I$ be a graded ideal of $K[x_1,\ldots,x_n]$ generated by homogeneous polynomials of a same degree $d$, and assume that $I$ has linear quotients. In this note, we use Horseshoe Lemma to give a relatively direct inductive construction of a minimal free resolution of $I$, which is called a $d$-linear resolution.}

\vs{3mm}\nin {\small Key Words:} {\small  Graded ideal; linear quotients; $d$-linear resolution; Horseshoe Lemma}

\vs{3mm}\nin {\small 2010 AMS Classification:} {\small   13F20, 05E40}

\end{minipage}
\end{center}

\section{Preliminaries}

\vs{4mm}We begin by recalling the following definition:

\begin{mydef} Let $I$ be a graded ideal of $S$. If there exists a sequence of homogeneous generators $f_1,
\ldots,f_m$ of $I$ such that for any $i$ with $1<i\leq m$, the colon ideal
$$\,\langle f_1,f_2,\ldots,f_{i-1}\,\rangle : f_i$$
 is generated by linear forms, then $I$ is said to have {\bf linear quotients} (with respect to the ordering $f_1,
\ldots,f_m$).
\end{mydef}\index{Graded ideals with linear quotients}

Note that {\it if a graded ideal $I$ is generated by linear forms, then
$I$ has linear quotients}. In fact, let $\,l_1,\ldots,l_r \,$ be a minimal generating set of linear forms of the ideal $I$.  We claim that {\it $\,l_1,\ldots,l_r$ is a regular sequence of $S$}. In fact, the minimality  ensures that it is an $\mathbb K$-independent subset of $S_1$.  Then $\,l_1,\ldots,l_r \,$ can be extended to a $\mathbb K$-basis, say $\,l_1,\ldots,l_n \,$, of $S_1$.
Then $x_i\mapsto    l_i,\, (\forall \, i=1,\ldots,n)$ induces an automorphism of the $\mathbb K$-algebra $S$. Since $\,x_1,\ldots,x_n \,$ is a regular sequence of $S$,  clearly $\,l_1,\ldots,l_n \,$ is also a regular sequence of $S$. In particular, $\,l_1,\ldots,l_r \,$ is a regular sequence, hence $$\,\langle l_1,l_2,\ldots,l_{i-1}\,\rangle : l_i=\langle l_1,l_2,\ldots,l_{i-1}\,\rangle, \,\forall 2\le i\le r.$$
Thus $I$ has linear quotients. In the case, note that $r=ht(I)$.

Recall that a polymatroidal  monomial ideal has linear quotients (see, e.g., \cite[Lemma 1.3]{HerzogTakayama}). For a simplicial complex $\Delta$, recall that the Stanley-Reisner ideal $I_{\Delta}$ has linear quotients iff  the Alexander dual complex  $\Delta^\vee$ of $\Delta$ is shellable (\cite[Proposition 8.2.5]{Herzog and Hibi}).

Now assume that a graded ideal $I$ has  linear quotients with respect to the ordering $f_1,
\ldots,f_m$. As in \cite{Herzog and Hibi 2006}, for each $2\le i\le m$, let $q_i(I)$ be the cardinal number of a minimal generating set of linear forms of the colon ideal $$\,\langle f_1,f_2,\ldots,f_{i-1}\,\rangle : f_i.$$ Then it follows from  the previous paragraph that $q_i(I)$ is independent of the choice of a minimal generating set of linear forms of the colon ideal. Let
$$q(I)=max\{q_i(I)\mid 2\le i\le m\}.$$
Note that if $I$ is generated by linear forms, then $q_i(I)=i-1$, thus $q(I)=r-1$, where $r$ is the cardinal number of a minimal generating set of linear forms of $I$.

For a graded ideal $I$ of the graded ring $S$, let
$$0\longra \oplus_{i=1}^{b_r} R(-d_{ri})\longra\cdots\longra  \oplus_{i=1}^{b_1} R(-d_{1i})\longra I\longra 0$$
be the minimal graded resolution of $I$  by free modules. The ideal $I$ is said to have a {\it pure resolution} if there are constants $d_1<d_2<\cdots<d_r$, such that
$$d_{1i}=d_1,\,\ldots,\, d_{ri}=d_r,\,\forall i.$$
If further $d_1=d, d_i=d_1+i-1,\,\forall 2\le i\le r$, then $I$ is said to {\it admit a $d$-linear resolution}. See \cite[Theorem 4.3.1]{Bruns and Herzog} for a characterization given by Eisenbud and Goto. For a simplicial complex $\Delta$, recall that the ideal $I_{\Delta}$ admits a linear resolution if and only if the complex $\Delta^\vee$ is Cohen-Macaulay over any field $K$ (Eagon-Reiner Theorem, see \cite[Theorem 8.1.9]{Herzog and Hibi}).

\section{A construction of minimal free resolution by using Horseshoe Lemma}

\begin{mylem} \label{LinearQuotientsInducesLinearResolution1} Let $I$ be a graded ideal of $S$ generated by linear forms. Then $I$ has $1$-linear resolution and $pd_S(I)=q(I)$.
\end{mylem}

\nin{\bf Proof.} Let $\,f_1,\ldots,f_r \,$ be a minimal generating set of linear forms of the ideal $I$ and denote $I_t=\langle \, f_1,\ldots,f_t\,\rangle, 1\le t\le r$. Then $q_t(I)=q(I_t)=t-1$ .

We prove the result by using induction on $r$. For $r=1$, the result is clear. Now assume $r>1$ and consider the exact sequence of graded $S$-modules and graded $S$-module homomorphisms
$$0\longra I_{r-1}\longra I_{r}\longra I_r/I_{r-1}\longra 0.$$
By induction, assume that $I_{r-1}$ has the following 1-linear free resolution of projective dimension $r-2$:
$$0\longra F_{r-2}\stackrel{\si_{r-2}}{\longrightarrow}\cdots\longra F_0\stackrel{\si_0}{\longrightarrow} I_{r-1}\longra 0.$$
Note that $I_r/I_{r-1}=S\ol{f_k}$, thus we have the following graded exact sequence
$$0\longra I_{r-1}(-1)\longra S(-1)\longra I_r/I_{r-1}\longra 0,$$
hence the following is the  1-linear free resolution of $I_r/I_{r-1}$ of projective dimension $r-1$:
$$0\longra F_{r-2}(-1)\stackrel{\si_{r-2}}{\longrightarrow}\cdots\longra F_0(-1)\stackrel{\si_0}{\longrightarrow}S(-1)\stackrel{\si_0'}{\longrightarrow} I_r/I_{r-1}\longra 0.$$

We then use the construction in proving the Horseshoe Lemma in homological algebra (see,e.g., \cite[Proposition 6.24]{Rotman}),  to construct a free resolution of $I_r$ with $pd_S(I_r)=r-1$:
\[\begin{array}{ccccccccccccc}
&&&& &0&&0&&\\
&&&& &\Big\downarrow & &\Big\downarrow\\
& 0 & \longra & 0 & \stackrel{j}{\longra} & 0\oplus F_{r-2}(-1) & \stackrel{\pi}{\longra} & F_{r-2}(-1) & \longra & 0 \\
&&&\vcenter{\llap{$$}}\,\Big\downarrow& &\Big\downarrow\vcenter{\rlap{$\delta_{r-1}$}} & &\Big\downarrow\vcenter{\rlap{$\si_{r-2}$}}\\
& 0 & \longra &F_{r-2} & \stackrel{j}{\longra} & F_{r-2}\oplus F_{r-3}(-1) & \stackrel{\pi}{\longra} & F_{r-3}(-1) & \longra & 0 \\
&&&\vcenter{\llap{$\si_{r-2}$}}\,\Big\downarrow& &\Big\downarrow\vcenter{\rlap{$\delta_{r-2}$}} & &\Big\downarrow\vcenter{\rlap{$\si_{r-3}$}}\\
&&&\vdots& &\vdots & &\vdots\\
&&&\Big\downarrow& &\Big\downarrow & &\Big\downarrow\\
& 0 & \longra & F_1 & \stackrel{j}{\longra} & F_1\oplus F_0(-1) & \stackrel{\pi}{\longra} & F_0(-1) & \longra & 0 \\
&&&\vcenter{\llap{$\si_1$}}\,\Big\downarrow& &\Big\downarrow\vcenter{\rlap{$\delta_1$}} & &\Big\downarrow\vcenter{\rlap{$\si_0$}}\\
& 0 & \longra & F_0 & \stackrel{j}{\longra} & F_0\oplus S(-1) & \stackrel{\pi}{\longra} & S(-1) & \longra & 0 \\
&&&\vcenter{\llap{$\si_0$}}\,\Big\downarrow& &\Big\downarrow\vcenter{\rlap{$\delta_0$}} & &\Big\downarrow\vcenter{\rlap{$\si_0'$}}\\
&0&\longra& I_{r-1}&\longra &I_{r}&\longra& I_r/I_{r-1}&\longra& 0\\
&&&\Big\downarrow& &\Big\downarrow & &\Big\downarrow\\
&&&0& &0&&0&&\\
\end{array}\]\par

By the construction, each $\delta_i$ is essentially the sum of $\si_i$ and $\si_i'$, thus the constructed middle resolution is a $1$-linear resolution of $I_r$. This shows that  $I$ has $1$-linear resolution and $pd_S(I)=r-1=q(I)$\qed

\vs{3mm}Again let $I_t=\langle \, f_1,\ldots,f_t\,\rangle, 1\le t\le r$. Let $L_k\,=\,\langle \,f_1,\ldots,f_{k-1}\,\rangle:f_k$ be the series of colon ideals. By Lemma \ref{LinearQuotientsInducesLinearResolution1}, $$pd_S(S/L_k)=pd_S(L_k)+1=q_k(I).$$
Note that
$$0\longrightarrow I_{k-1}\longrightarrow I_k\longrightarrow(S/L_k)(-d)\longrightarrow 0$$
is a exact sequence of graded modules, thus $$pd(I_k)=max\{pd(I_{k-1}),\, q_k(I)\}=q(I_k).$$
So, if we use the previous lemma and the proof to the above sequences, then in a similar manner we use induction and the Horseshoe Lemma to give a relatively direct proof to the following:

\begin{mthm} (\cite[{\textrm Proposition} $8.2.1$  {\textrm and}  {\textrm Corollary} 8.2.2]{Herzog and Hibi}) Let $I$ be  a graded ideal  of $S$ generated by homogeneous polynomials of   degree $d$.  If $I$ has linear quotients, then $I$ has a $d$-linear resolution and, $pd_S(I)=q(I)$.
\end{mthm}\index{monomial ideal!linear quotients}\index{monomial ideal!linear resolution}


\begin{thebibliography}{gg}


\bibitem{Bruns and Herzog} Bruns W. and Herzog J.. {\it Cohen-Macaulay Rings.} {\bf Cambrisge University Press}, Cambridge, Rev. Ed., 1997.

\bibitem{Eisenbud}  Eisenbud D. {\it Commutative Algebra with a View Toward Algebraic Geometry}. {\bf Springer Science + Business Media, Inc.}
 2004.

\bibitem{Herzog and Hibi 2006} Herzog J.  and  Hibi T. Cohen-Macaulay polymatroidal ideals, European Journal of Combinatorics $27 (2006) \,\,513-517.$

\bibitem{Herzog and Hibi} Herzog J.  and  Hibi T.  {\it Monomial Ideals}.  {\bf Springer-Verlag London Limited}, $2011$.

\bibitem{HerzogTakayama} J. Herzog, Y. Takayama£¬ Resolutions by mapping cones, Homology Homotopy Appl. $4 (2002)\,\, 277-294.$

\bibitem{Rotman} Rotman J.J. {\it An Introduction To Homological Algebra}.  {\bf Springer Science +Business Media}, LLC 2009. Second Edition.


\end{thebibliography}
\end{document}